\documentclass[preprint]{imsart}

\RequirePackage[OT1]{fontenc}
\RequirePackage{amsthm,amsmath}
\RequirePackage[numbers]{natbib}
\RequirePackage[colorlinks,citecolor=blue,urlcolor=blue]{hyperref}

\startlocaldefs
\numberwithin{equation}{section}
\theoremstyle{plain}

\newtheorem{defn}{Definition}[section]

\endlocaldefs

\begin{document}

\begin{frontmatter}
\title{Satisfaction Problem  of Consumers  Demands measured by ordinary  ``Lebesgue
measures'' in $R^{\infty}$}
\runtitle{Satisfaction Problem  of Consumers  Demands}

\begin{aug}
\author{\fnms{Gogi} \snm{Pantsulaia}\thanksref{t1}\ead[label=e1]{g.pantsulaia@gtu.ge}}
\and
\author{\fnms{Givi } \snm{Giorgadze}\ead[label=e2]{g.giorgadze@gtu.ge}}
\affiliation{}

\thankstext{t1}{The present work  was partially supported on the Shota Rustaveli
National Science Foundation Grants: $\#$ GNFS 31/25,~$\#$ GNFS / FR  1165-100 / 14.
}
\runauthor{G.Pantsulaia and G.Giorgadze}

\address{I.Vekua Institute of Applied Mathematics, Tbilisi - 0143, Georgian Republic \\
\printead{e1}\\
Georgian Technical University, Tbilisi - 0175, Georgian Republic\\
  \printead{e2}
}
\end{aug}

\begin{abstract}
 In the present paper we consider the following Satisfaction Problem  of Consumers  Demands (SPCD): {\it The supplier  must supply the measurable system of the measure $m_k$  to  the $k$-th  consumer at time  $t_k$   for $1 \le k \le n$. The  measure   of the supplied measurable system is changed  under action of  some  dynamical system;
What is a minimal measure  of measurable system  which must take the supplier  at the initial time  $t=0$  to satisfy  demands of all consumers ? }
In this paper we consider Satisfaction Problem  of Consumers  Demands measured by ordinary ``Lebesgue
measures'' in $R^{\infty}$ for various  dynamical systems in $R^{\infty}$. In order to solve this problem  we  use  Liouville type theorems for them which  describes the dependence between initial  and resulting  measures of the entire system.
\end{abstract}

\begin{keyword}[class=MSC]
\kwd[Primary ]{37-xx}
\kwd{28Axx}
\kwd[; Secondary ]{28C10}
{28D10}
\end{keyword}

\begin{keyword}
\kwd{Foerster–-Lasota equation}
\kwd{Black - Scholes equation}
\kwd{ordinary  ``Lebesgue
measure''}
\end{keyword}

\end{frontmatter}

\section{Formulation of the main problem}

     In a wide class of general  systems  the so-called dynamical systems can be  separated 
which, as usual,   are used  for  description  of the behaviour of various  physical,
economic and social processes  during the time (see, for example \cite{NemStep49}). 

     Let $(E,\rho)$  be a metric space.  Recall, that a family of mappings $(\Phi_t)_{t \in R}$ with
$\Phi_t : E \to E$ for $t \in R$  is called a dynamical system  if
it satisfies the following three conditions:

 1)  $\Phi_0(x)=x$  for each element  $x \in E$ ;

 2)  A mapping $\Phi :  E \times R \to E$ defined by $\Phi(x,t)=\Phi_t(x)$
 is continuous with respect to the variables $x$ and $t$;

 3)  if $x \in E$, $t_1 \in R$   and  $t_2 \in R$,  then
 $\Phi_{t_1}(\Phi_{t_2}(x))=\Phi_{t_1+t_2}(x).$

\medskip

{\bf Satisfaction  Rule   of Consumers  Demands  under  Dynamical Sistem):}  Let $(\Phi_t)_{t \in R}$ be some dynamical system defined in a
metric space $(E,\rho)$  and let $\nu$  be  a Borel measure on
$E$.  It is assumed that $t_1 < t_2 < \cdots < t_n$  and  $m_k>0$ for $1 \le k \le m$.  Supplier  is a person who choose any Borel subset $Y \subseteq E$ of positive $\nu$  measure  called an initial system
and  satisfies $n$ consumers demands as follows:

(i)~ At time  $t_1$  supplier  gives   a set   $\Phi_{t_1}(Y)$ and   $1$-th consumer  at time  $t_1$    choose  a Borel subset  $C_1 \subseteq  \Phi_{t_1}(Y)$ of the measure $m_1$ (called $1$-th  demand)  for which  $\nu(C_1)=m_1$ as follows:
$$C_1= \tau\{ X : X \in B( \Phi_{t_1}(Y))~\&~\nu(X)=m_1\},$$
where  $B( \Phi_{t_1}) $  denotes the class of all Borel subsets  of $ \Phi_{t_1}(Y)$  and   $\tau$ denotes  an operator of global choice.

(ii)~ At time  $t_2$ supplier gives   a set   $\Phi_{t_2 -t_1}(\Phi_{t_1}( Y)  \setminus  C_1 )$ and   $2$-th consumer  at time $t_2$    choose  a Borel subset  $C_2 \subseteq  \Phi_{t_2 -t_1}(\Phi_{t_1}( Y)  \setminus  C_1 )$ of the measure $m_2$ (called $2$-th  demand)  for which  $\nu(C_2)=m_2$ as follows:
$$C_2= \tau\{ X : X \in B(\Phi_{t_2 -t_1}(\Phi_{t_1}( Y)  \setminus  C_1 ))~\&~\nu(X)=m_2\},$$

and so on.

\medskip

{\bf Satisfaction Problem  of Consumers  Demands ({\bf SPCD} ):}  {\it Assume that the supplier  must satisfy  $n$ consumers demands  by the rule described above.   What minimal measure  of a measurable system  which must take the supplier  at the initial time  $t=0$  for  satisfaction  demands of all consumers ? }

\medskip

We   plan to consider  {\bf SPCD}  in the case, when  $(E,\rho)$  is an infinite-dimensional topological vector space $R^{\infty}$ equipped with Tychonov metric,  $\nu$ is any ordinary ``Lebesgue
measure'' in $R^{\infty}$ (cf. \cite{Pan09-ord})  and   $(\Phi_t)_{t \in R}$ in $R^{\infty}$   is  a dynamical system  defined  by the one from the following  differential  equations:

$\bullet$~  von Foerster–-Lasota differential
equation in ${\bf R}^{\infty}$  (cf. \cite{Pan-Gio12(2)}) ;

$\bullet$~ The  Black–-Scholes equation (cf. \cite{Pan-Giorgadze2014});

$\bullet$~ Infinite generalised Maltusian growth equation in $R^{\infty}$ (cf. \cite{Pan07});

$\bullet$~ Fourier differential  equation (cf. \cite{Pan-Giorgadze2014}).

The rest of the paper is the following.

In Section 2 we give constructions of  ordinary ``Lebesgue
measures'' in $R^{\infty}$ . In the next sections we discuss the Satisfaction Problem  of Consumers  Demands
for above mentioned mathematical models.

\section{Auxiliary notions and propositions from measure
theory and linear algebra}

\medskip

Let $(\beta_j)_{j \in N} \in [0, +\infty]^N$.

\medskip

\noindent {\bf Definition 2.1}~ We say that a number $\beta \in
[0, +\infty]$ is an ordinary product of numbers $(\beta_j)_{j \in
N}$ if
$$
\beta=\lim_{n \to \infty}\prod_{i=1}^n\beta_i.
$$
An ordinary product of numbers $(\beta_j)_{j \in N}$ is
denoted by ${\bf (O)}\prod_{i \in N}\beta_i$.

\medskip

Let $\alpha=(n_k)_{k \in N} \in (N \setminus
\{0\})^{N}$. We set
$$
F_0=[0,n_0] \cap N ,~F_1=[n_0\!+\!1,n_0\!+\!n_1]\cap
N,~\dots, $$
$$
F_k=[n_0\!+\cdots+ n_{k-1}\!+\!1,n_0\!+\cdots
+\!n_k]\cap N, \dots\,.
$$
\medskip

\noindent {\bf Definition 2.2}~ We say that a number $\beta \in
[0, +\infty]$ is an ordinary $\alpha$-product of numbers
$(\beta_i)_{i \in N}$ if $\beta$ is   an ordinary product
of numbers $(\prod_{i  \in F_k}\beta_i)_{k \in N}$. An
ordinary $\alpha$-product of numbers $(\beta_i)_{i \in
N}$ is denoted by ${\bf (O,\alpha)}\prod_{i \in
N}\beta_i$.

\medskip

\noindent {\bf Definition 2.3}~ Let $\alpha=(n_k)_{k \in
N} \in (N \setminus \{0\})^{N}$. Let
$\mathcal{(\alpha)OR} $ be the class of all infinite-dimensional
measurable $\alpha$-rectangles  $R = \prod_{i \in N}R_i
(R_i \in \mathcal{B}({\bf R}^{n_i})$) for which an ordinary
product of numbers $(m^{n_i}(R_i))_{i \in N}$ exists and
is finite, where $m$ denotes a linear Lebesgue measure in $R$.

\medskip

\noindent {\bf Definition 2.4}~We say that a  measure  $\lambda$ being the completion of a
translation-invariant Borel measure  is an ordinary
$\alpha$-Lebesgue measure on $R^{\infty}$(or, shortly,
O$(\alpha)$LM)  if ~for every $R  \in \mathcal{(\alpha)OR} $ we
have
$$
\lambda(R)={\bf (O)}\prod_{k \in N}m^{n_k}(R_k).
$$

\medskip

\noindent {\bf Lemma 2.1}~(\cite{Pan07},  Theorem 1, p. 216 )~{\it
For every $\alpha=(n_i)_{i \in N} \in (N
\setminus \{0\})^{N}$, there exists a Borel measure
$\mu_{\alpha}$ on ${\bf R}^{\infty}$ which is O$(\alpha)$LM.}

\medskip

\noindent {\bf Lemma 2.2}~(\cite{Pan09}, Theorem 3, p. 9.)~{\it
Let $\alpha=(n_i)_{i \in N} \in (N \setminus
\{0\})^{N}$, and let $T^{n_i} : {\bf R}^{n_i} \to
{\bf R}^{n_i}, i
> 1$, be a  family of linear transformation with Jacobians
$\Delta_i \neq 0$ and $0 < \prod_{i=1}^{\infty}\Delta_i < \infty$.
Let $T^{N} : R^{N} \to R^{N}$ be the map defined
by
$$T^{N}(x) =
(T^{n_1}(x_1, \cdots, x_{n_1}),T^{n_2}(x_{n_1+1}, \cdots,
x_{n_1+n_2}), \cdots),$$ where $x = (x_i)_{i \in N} \in
R^{N}$. Then for each $ E \in
\mathcal{B}(R^{N}),$ ~we have
$$\mu_{\alpha}(T^{N}(E)) =
(\prod_{i=1}^{\infty}\Delta_i) \mu_{\alpha}(E).$$}

\medskip

In context with another interesting properties of partial analogs
of the Lebesgue measures in ${\bf R}^{\infty}$, the reader can
see \cite{Bak91}, \cite{Bak04}, \cite{Pan04-1},
\cite{Pan09-ord}, \cite{Pan10},\cite{Pan07}.

\medskip

In the sequel we identify the vector space $R^{\infty}$  of all
real-valued sequences with  the vector space of all real-valued
infinite-dimensional vector-columns.

The need the following auxiliary proposition
from linear algebra.

\medskip

\noindent {\bf Lemma 2.3}~( \cite{Gant66}, \S 6, Section 1)~{\it
Let $\alpha=(n_i)_{i \in N} \in (N \setminus
\{0\})^{N}$ and, let $A=(a_{ij})_{i,j \in N}$ be an
infinite-dimensional real-valued $\alpha$-cellular matrix. Let us
consider a linear autonomous differential equation of the first
order
$$
\frac{d}{d t} ((a_k)_{k \in N})=A \times (a_k)_{k \in N} \eqno
(2.1)
$$
with an initial condition
$$(a_k(0))_{k \in N}=(c_k)_{k \in N}\in {\bf R}^{\infty}. \eqno (2.2)
$$
Then the solution of $(2.1)-(2.2)$ is given by
$$
(a_k(t))_{k \in N}=\exp(tA)\times (c_k)_{k \in N}. \eqno (2.3) $$}

\noindent \begin{proof} Let us present the column $(a_k(t))_{k \in
N}$ in the Maclaurin series as follows:
$$
(a_k(t))_{k \in N}=\sum_{m=0}^{\infty}\frac{(a^{(m)}_k(0))_{k \in
N}}{m!}t^m. \eqno (2.4)
$$

Take  into account the validity of the formula
$$
(a^{(m)}_k(0))_{k \in N}=(\frac{d^m a_k(t)}{d t^m}|_{t=0})_{k \in
N}=A^m \times (a_k(0))_{k \in N},~~(2.5)
$$
we get
$$
(a_k(t))_{k \in N}=\sum_{m=0}^{\infty}\frac{(a^{(m)}_k(0)t^m)_{k
\in N}}{m!}= \sum_{m=0}^{\infty}\frac{(tA)^m}{m!} \times
(a_k(0))_{k \in N}=\exp(tA)\times (c_k)_{k \in N} ~~\eqno (2.6)
$$
\end{proof}

\medskip

In the sequel  we will need  some notions characterizing the
behavior of some dynamical systems $(\Phi_t)_{t \ge 0}$ in ${\bf R}^{\infty}$.

Let  $\nu$ be any `` Lebesgue measure'' in ${\bf R}^{\infty}$
(see, for example, \cite{Pan04-1}, \cite{Pan09-ord}).

\medskip

\noindent {\bf Definition 2.5}~ We say that the dynamical system
$(\Phi_t)_{t \ge 0}$ is stable in the sense of a `` Lebesgue
measure'' $\nu$  if  the flow  preserves the measure $\nu$, i.e.
$$
(\forall t)(\forall D)(0 < t < \infty ~\&~0< \nu(D)< \infty
\rightarrow \nu(\Phi_t(D))=\nu(D)).~~\eqno (2.7)
$$

\medskip

\noindent {\bf Definition  2.6} We say that the dynamical system
$(\Phi_t)_{t \ge 0}$ is expansible in the sense of a ``Lebesgue
measure'' $\nu$  if
$$
(\forall t_1)(\forall t_2)(\forall D)(0 < t_1<t_2< \infty ~\&~0<
\nu(D)< \infty \rightarrow \nu(\Phi_{t_1}(D))<
\nu(\Phi_{t_2}(D))).~~\eqno (2.8)
$$
\medskip

\noindent {\bf Definition  2.7} We say that the dynamical system
$(\Phi_t)_{t \ge 0}$ is pressing in the sense of a ``Lebesgue
measure'' $\nu$  if  the 
flow is dissipative, i.e., 
$$
(\forall t_1)(\forall t_2)(\forall D)(0 < t_1<t_2< \infty ~\&~0<
\nu(D)< \infty \rightarrow \nu(\Phi_{t_1}(D))>
\nu(\Phi_{t_2}(D))).~~\eqno (2.9)
$$
\medskip

\noindent {\bf Definition  2.8} We say that the dynamical system
$(\Phi_t)_{t \ge 0}$ is totally expansible in the sense of a
``Lebesgue measure'' $\nu$  if
$$
(\forall t)(\forall D)(0 < t< \infty ~\&~0< \nu(D)< \infty
\rightarrow \nu(\Phi_{t}(D))=+\infty ).~~\eqno (2.10)
$$
\medskip

\noindent {\bf Definition  2.9} We say that the dynamical system
$(\Phi_t)_{t \ge 0}$ is totally pressing in the sense of a
``Lebesgue measure'' $\nu$  if
$$
(\forall t)(\forall D)(0 < t< \infty ~\&~0< \nu(D)< \infty
\rightarrow \nu(\Phi_{t}(D))=0).~~\eqno (2.11)
$$
\medskip

\section{ Satisfaction Problem  of Consumers  Demands in von Foerster–-Lasota model in ${\bf R}^{\infty}$}

In this section we consider  a certain concept \cite{Pan-Gio13(1)}
for a solution of some differential equations by `` Maclaurin
Differential Operators'' in ${\bf R}^{\infty}$.

{\bf Definition 3.1} `` Maclaurin  differential operator''
$(\mathcal{M})\frac{\partial}{\partial x}$ in ${\bf R}^{\infty}$
is defined as follows:
$$
(\mathcal{M}) \frac{\partial}{\partial x}(\begin{pmatrix}
   a_0&  \\
   a_1&  \\
   a_2&  \\
  \vdots &
\end{pmatrix})=\begin{pmatrix}
   0&  1&  0&  0&   \ldots & \\
   0&  0&  2&  0&   \ddots & \\
   0&  0&  0&  3&  \ddots &  \\
   0&  0&  0&  0&  \ddots &  \\
    \vdots &   \vdots &   \vdots &   \vdots &  \ddots &
\end{pmatrix} \times \begin{pmatrix}
  a_0 \\
   a_1\\
   a_2\\
  \vdots
\end{pmatrix}.\eqno(3.1)
$$

{\bf Definition 3.2} `` Maclaurin  differential operator''
$(\mathcal{M}){x \frac{\partial}{\partial x}}$  in ${\bf
R}^{\infty}$ is defined as follows:
$$
(\mathcal{M}) x \frac{\partial}{\partial x} ((a_k)_{k \in N})=
\begin{pmatrix}
   0&  0&  0&  0&   \ldots&  \\
   0&  1&  0&  0&   \ddots&  \\
   0&  0&  2&  0&  \ddots&  \\
   0&  0&  0&  3&  \ddots&  \\
    \vdots&   \vdots&   \vdots&   \ddots&  \ddots&
\end{pmatrix} \times \begin{pmatrix}
  a_0 \\
   a_1\\
   a_2\\
  \vdots
\end{pmatrix}. \eqno(3.2)
$$

 {\bf Definition 3.3} `` Maclaurin  differential operator'' $(\mathcal{M})x^2
\frac{\partial^2}{\partial x^2}$  in ${\bf R}^{\infty}$ is defined
as follows:
$$
(\mathcal{M}){ x^2 \frac{\partial^2}{\partial x^2}} ((a_k)_{k \in
N})=
\begin{pmatrix}
    0&   0&   0&   0&    \ldots&   \\
    0&   0&   0&   0&    \ddots&   \\
    0&   0&   1\times 2&   0&   \ddots&   \\
    0&   0&   0&   2\times 3&   \ddots&   \\
     \vdots&    \vdots&    \vdots&    \ddots&   \ddots&
\end{pmatrix} \times \begin{pmatrix}
   a_0 \\
    a_1\\
    a_2\\
   \vdots
\end{pmatrix}. \eqno(3.3)
$$

 {\bf Definition 3.4}  Formally,  we set that the factorial of each negative integer number is equal to $+\infty$.  Then `` Maclaurin  differential operator'' $(\mathcal{M}){x^n
\frac{\partial^n}{\partial x^n}}$  in ${\bf R}^{\infty}$    is defined as
follows:
$$
(\mathcal{M}){ x^n \frac{\partial^n}{\partial x^n}} ((a_k)_{k \in
N})= \begin{pmatrix}
    \frac{0!}{(0-n)!}&   0&   0&   0&    \ldots&   \\
    0&   \frac{1!}{(1-n)!}&   0&   0&    \ddots&   \\
    0&   0&   \frac{2!}{(2-n)!}&   0&   \ddots&   \\
    0&   0&   0&   \frac{3!}{(3-n)!}&   \ddots&   \\
     \vdots&    \vdots&    \vdots&    \ddots&   \ddots&
\end{pmatrix} \times \begin{pmatrix}
   a_0 \\
    a_1\\
    a_2\\
   \vdots
\end{pmatrix}. ~~\eqno(3.4)
$$

 {\bf Theorem 3.3}~(\cite{Pan(14)},Theorem 11.2, p.139)~ {\it  Let $(A_n)_{0 \le n \le m}(m \in N)$ be
a sequence of real numbers.
 Let consider a non-homogeneous  `` Maclaurin  differential operators'' equation  of the first order
$$
{\frac{d}{d t}}((a_k)_{k \in N})=  \sum_{n=0}^m A_n (\mathcal{M})
x^n \frac{\partial^n}{\partial x^n} \times (a_k)_{k \in N}+
(f_k(t))_{k \in N}     \eqno(3.5)
$$
with initial condition
$$
(a_k(0))_{k \in N}=(C_k)_{k \in N},  \eqno(3.6)
$$
where

(i) $(C_k)_{k \in N} \in {\bf R}^{\infty}$;

(ii) $(f_k(t))_{k \in N}$ is the sequence of continuous functions
on ${\bf R}$.

Then $$(a_k(t))_{k \in N}=\big(e^{t\sum_{n=0}^m
A_n\frac{k!}{(k-n)!}}C_k+ \int_{0}^te^{(t-\tau)\sum_{n=0}^m
A_n\frac{k!}{(k-n)!}}f_k(\tau)d\tau \big)_{k \in N}.\eqno(3.7)
$$}

\medskip

We have the following consequence of Theorem 3.3.

{\bf Corollary  3.1 }~~{\it Let consider the von Foerster–-Lasota operator
equation in ${\bf R}^{\infty}$ defined by
$$
{\frac{d}{d t}}((a_k)_{k \in N})= -(\mathcal{M})\big( x \frac{\partial}{\partial x}\big)((a_k)_{k \in N})+ \gamma (a_k)_{k \in N} \eqno(3.8)
$$
with initial condition
$$
(a_k(0))_{k \in N}=(C_k)_{k \in N}  \in {\bf R}^{\infty}. \eqno(3.9)
$$

Then $$(a_k(t))_{k \in N}=\big(e^{t(\gamma - k)} C_k \big)_{k \in N}.\eqno(3.10)
$$
}

{\bf Satisfaction Problem  of Consumers  Demands in von Foerster–-Lasota model.}

\medskip

  Let consider $(1,1, \cdots)$-ordinary
''Lebesgue measure''  $\mu_{(1,1, \cdots)}$. By Lemma 2.2 we know that $\mu_{(1,1, \cdots)}(\Phi_t(X))= e^{\sum_{k =0}^{\infty}t(\gamma - k)}\mu_{(1,1, \cdots)}(X)$, where the von Foerster–-Lasota motion $\Phi_t : R^{\infty} \to R^{\infty}$ is defined by

$$\Phi_t \big( (C_k)_{k \in N}\big)= \big(e^{t(\gamma - k)} C_k \big)_{k \in N}\eqno(3.11) $$
for $(C_k)_{k \in N} \in  R^{\infty}$.

Since $e^{\sum_{k =0}^{\infty}t(\gamma - k)}=0$, we claim for an arbitrary initial system $S_0 \in \mathcal{B}(R^{\infty})$ and $t >0$, the set
$\mu_{(1,1, \cdots)}(\Phi_t(S_0))=0$. Hence the first consumer  can not choose a Borel subset $C_1 \subseteq S_{t_1}= \Phi_{t_1}(S_0)$ for which $\mu_{(1,1, \cdots)}(C_1)=m_1>0$. The latter relation means that Satisfaction Problem  of Consumers  Demands in von Foerster–-Lasota model has no any solution.

\medskip

\section{Satisfaction Problem  of Consumers  Demands in  Black –- Scholes  Model  in $R^{\infty}$ }

The Black–-Scholes differential  equation  in $R^{\infty}$  has the following form :

$${\frac{d}{d t}}(a_k)_{k \in N}= - \frac{1}{2}\sigma^2 (\mathcal{M})\big( x^2 \frac{\partial}{\partial x^2}\big)((a_k)_{k \in N})
-r(\mathcal{M})\big( x \frac{\partial}{\partial x}\big)((a_k)_{k \in N})+r (a_k)_{k \in N}\eqno (4.1)
$$

Notice that $(4.1)$ is a particular case of $(3.5)$ for which
$m=2$,  $A_0=r$, $A_1=-r$, $A_2=- \frac{1}{2}\sigma^2$. Following
Theorem 3.3, the solution of $(4.1)$ has the form
$$
(a_k(t))_{k \in N}=\big(e^{t(r\frac{k!}{(k-0)!}-r
\frac{k!}{(k-1)!}-
 \frac{1}{2}\sigma^2\frac{k!}{(k-2)!})}C_k\big)_{k \in N}\eqno(4.2)
 $$

 {\bf Satisfaction Problem  of Consumers  Demands in Black–-Scholes model.}

 \medskip

  Let consider $(1,1, \cdots)$-ordinary
''Lebesgue measure''  $\mu_{(1,1, \cdots)}$. By Lemma 2.2 we know that $\mu_{(1,1, \cdots)}(\Phi_t(X))= e^{\sum_{k =0}^{\infty}t(\gamma - k)}\mu_{(1,1, \cdots)}(X)$, where the von Foerster–-Lasota motion $\Phi_t : R^{\infty} \to R^{\infty}$ is defined by

$$\Phi_t \big( (C_k)_{k \in N}\big)= \big(e^{t(r\frac{k!}{(k-0)!}-r
\frac{k!}{(k-1)!}-
 \frac{1}{2}\sigma^2\frac{k!}{(k-2)!})}C_k\big)_{k \in N}\eqno(4.3) $$
for $(C_k)_{k \in N} \in  R^{\infty}$.

Since $\sum_{k \in N}r\frac{k!}{(k-0)!}-r
\frac{k!}{(k-1)!}-
 \frac{1}{2}\sigma^2\frac{k!}{(k-2)!}=-\infty,$ we claim for an arbitrary initial system $S_0 \in \mathcal{B}(R^{\infty})$ and $t >0$, the set
$\mu_{(1,1, \cdots)}(\Phi_t(S_0))=0$. Hence first consumer  can not choose a Borel subset $C_1 \subseteq S_{t_1}= \Phi_{t_1}(S_0)$ for which $\mu_{(1,1, \cdots)}(C_1)=m_1>0$. The latter relation means that Satisfaction Problem  of Consumers  Demands in Black–-Scholes model has no any solution.

 \section{Satisfaction Problem  of Consumers  Demands in  infinite continuous generalised Malthusian growth model  in $R^{\infty}$ }

\medskip

Let us consider an infinite non-antagonistic family of populations and let $\Psi_k(t)$ be the population function of the $k$-th Population. Then the generalised continuous Malthusian growth model for an infinite family of non-antagonistic populations is described by the following linear differential equation
$$
\frac{d ((a_k(t))_{k \in N}}{dt}=A \times (a_k(t))_{k \in N} \eqno(5.1)
$$
with an initial condition $$(a_k(0))_{k \in N}=(a_k)_{k \in N} \in R^{\infty}, \eqno(5.2) $$ where
$A$ is an infinite-dimensional real-valued diagonal matrix with diagonal elements $(\lambda_k)_{k \in N}$.

By  Lemma 2.3 we know that the solution of (5.1) is given by

$$(a_k(t))_{k \in N}=(e^{t\lambda_k}a_k)_{k \in N}.\eqno(5.3)$$

{\bf Satisfaction Problem  of Consumers  Demands in infinite continuous generalised Malthusian growth model  in $R^{\infty}$ .}

 \medskip

 Let consider $(1,1, \cdots)$-ordinary
''Lebesgue measure''  $\mu_{(1,1, \cdots)}$. By Lemma 2.2 we know that $$\mu_{(1,1, \cdots)}(\Phi_t(X))= e^{t\sum_{k =0}^{\infty}\lambda_k}\mu_{(1,1, \cdots)}(X),\eqno(5.4)$$ where the infinite continuous generalised Malthusian growth motion $\Phi_t : R^{\infty} \to R^{\infty}$ is defined by

$$\Phi_t \big( (a_k)_{k \in N}\big)= (e^{t\lambda_k}a_k)_{k \in N}\eqno(5.5) $$
for $(a_k)_{k \in N} \in  R^{\infty}$.

{\bf Case 1.} $\sum_{k=0}^{\infty}\lambda_k$ is divergent. In that case we do not know whether SPCD in infinite continuous generalised Malthusian growth model has any solution.

{\bf Case 2.} $\sum_{k=0}^{\infty}\lambda_k=-\infty$. In that case  SPCD in infinite continuous generalised Malthusian growth model has no any solution.

{\bf Case 3.} $\sum_{k=0}^{\infty}\lambda_k=+\infty$. In that case SPCD in infinite continuous generalised Malthusian growth model also has no any solution because if the supplier take an arbitrary
measurable system  of the positive $(1,1, \cdots)$-ordinary
''Lebesgue measure''  $\mu_{(1,1, \cdots)}$, then  demands of all consumers
will be always  satisfied.

{\bf Case 4.}$\sum_{k=0}^{\infty}\lambda_k$ is convergent and $-\infty < \sum_{k=0}^{\infty}\lambda_k< +\infty$.

Let show that the solution of  SPCD in infinite continuous generalised Malthusian growth model is defined by
$$
m=\sum_{k=1}^n m_ke^{-t_k\sum_{i=1}^{\infty}\lambda_k}\eqno(5.6).
$$

Indeed, let choose an arbitrary Borel subset $S_0 \subset R^{\infty}$ with $\mu_{(1,1, \cdots)}(S_0)=m$. At the moment $t=t_1$ the set $S_0$
is transformed into set $\Phi_{t_1}(S_0)$ whose $\mu_{(1,1, \cdots)}$ measure is equal to  $$e^{t_1\sum_{i=1}^{\infty}\lambda_k}m=e^{t_1\sum_{i=1}^{\infty}\lambda_k} (\sum_{k=1}^n m_k e^{-t_k\sum_{i=1}^{\infty}\lambda_k})=$$
$$
\sum_{k=1}^n m_k e^{(t_1-t_k)\sum_{i=1}^{\infty}\lambda_k}=m_1 +\sum_{k=2}^n m_k e^{(t_1-t_k)\sum_{i=1}^{\infty}\lambda_k}\eqno(5.7).$$
When the demand $C_1$ of the first consumer will be satisfied, we obtain the set $\Phi_{t_1}(S_0) \setminus C_1$ for which
  $$\mu_{(1,1, \cdots)}(\Phi_{t_1}(S_0) \setminus C_1)= \sum_{k=2}^n m_k e^{(t_1-t_k)\sum_{i=1}^{\infty}\lambda_k}\eqno(5.8).$$
At moment $t=t_2$  the set $\Phi_{t_1}(S_0) \setminus C_1$ is transformed into set  $\Phi_{t_2-t_1}(\Phi_{t_1}(S_0) \setminus C_1)$ whose  $\mu_{(1,1, \cdots)}$  measure is equal to $$e^{(t_2-t_1)\sum_{i=1}^{\infty}\lambda_k}\mu_{(1,1, \cdots)}(\Phi_{t_1}(S_0) \setminus C_1)=
$$
$$
e^{(t_2-t_1)\sum_{i=1}^{\infty}\lambda_k}\sum_{k=2}^n m_k e^{(t_1-t_k)\sum_{i=1}^{\infty}\lambda_k}=
\sum_{k=2}^n m_k e^{(t_2-t_1)+(t_1-t_k)\sum_{i=1}^{\infty}\lambda_k}=
$$
$$
\sum_{k=2}^n m_k e^{(t_2-t_k)\sum_{i=1}^{\infty}\lambda_k}=
m_2+\sum_{k=3}^n m_k e^{(t_2-t_k)\sum_{i=1}^{\infty}\lambda_k}\eqno(5.9).
$$
When the demand $C_2$ of the second consumer will be satisfied, we obtain the set $\Phi_{t_2-t_1}(\Phi_{t_1}(S_0) \setminus C_1)\setminus C_2$ for which $$\mu_{(1,1, \cdots)}(\Phi_{t_2-t_1}(\Phi_{t_1}(S_0) \setminus C_1)\setminus C_2)= \sum_{k=3}^n m_k e^{(t_2-t_k)\sum_{i=1}^{\infty}\lambda_k}, \eqno(5.10)$$
and so on.

Now it obvious that at the moment $t=t_n$ we obtain a set whose $\mu_{(1,1, \cdots)}$ measure  exactly coincides with the positive number $m_n$ and hence the demand of the $n$-th consumer will be satisfied.

{\bf Observation 5.1} We have showed that Satisfaction Problem  of Consumers  Demands in infinite continuous generalised Malthusian growth model  in $R^{\infty}$ has the solution if dynamical system defined by (5.3) is  pressing, expansible or stable  in the sense of the measure  $\mu_{(1,1, \cdots)}$.  When  $(\Phi_t)_{t \in R}$ is  totally pressing or totally expansible the the same problem has no any solution.

\medskip

\section{Satisfaction Problem  of Consumers  Demands for dynamical system  defined  by  the Fourier differential equation
in $R^{\infty}$ }

\begin{defn}

''Fourier   differential operator'' $(\mathcal{F}){\frac{d}{d x}}
: R^{\infty} \to R^{\infty} $  is defined as follows:

$$(\mathcal{F}){\frac{d}{d x}}\big(
\begin{pmatrix}
  \frac{a_0}{2} \\
   a_1\\
   b_1\\
   a_2\\
   b_2\\
   a_3\\
   b_3\\
  \vdots
\end{pmatrix}\big)=\begin{pmatrix}
   0&  0&  0& 0 & 0 &  0&  0& \dots&  \\
   0&  0&  \frac{1\pi}{l}&  0&  0&  0&  0& \dots&  \\
   0&  -\frac{1\pi}{l}&  0&  0&  0&  0&  0& \dots&  \\
   0&  0&  0&  0&  \frac{2\pi}{l}& 0 &  0& \dots&  \\
   0& 0 &  0&  -\frac{2\pi}{l}&  0&  0&  0&  \dots& \\
   0&  0&  0&  0&  0&  0&  \frac{3\pi}{l}& \ddots&  \\
   0&  0&  0&  0&  0&  -\frac{3\pi}{l}&  0& \ddots&  \\
   \vdots&  \vdots&  \vdots&  \vdots&  \vdots&  \vdots&  \ddots& \ddots&
\end{pmatrix}   \times \begin{pmatrix}
  \frac{a_0}{2} \\
   a_1\\
   b_1\\
   a_2\\
   b_2\\
   a_3\\
   b_2\\
  \vdots
\end{pmatrix}. \eqno(6.1)
$$
\end{defn}

Suppose  that $(A_n)_{n \in N} \in R^{N}$  be  a sequence of real numbers such that
$$\sigma_k=\sum_{n=0}^{\infty} (-1)^n A_{2n}(\frac{k \pi}{l})^{2n} \eqno(6.2)$$
 and
 $$ \omega_k=\sum_{n=0}^{\infty} (-1)^{n}A_{2n+1}(\frac{k\pi}{l})^{2n+1}\eqno(6.3)$$
are convergent  for each $k \ge 1$.

{\bf Corollary 6.1.}(cf. \cite{Pan2012}) ~Let us consider a partial differential equation of the
first order
$$
{\frac{\partial}{\partial t}}((a_k)_{k \in N})=
\Big(\sum_{n=0}^{\infty}
A_n\Big((\mathcal{F})\frac{\partial}{\partial x}\Big)^{n}\Big)
\times (a_k)_{k \in N} \eqno(6.4)
$$
with initial condition
$$
(a_k(0))_{k \in N}=(C_k)_{k \in N}, \eqno(6.5)
$$
where $(C_k)_{k \in N} \in {\bf R}^{\infty}$;

Suppose that the sequence of real numbers $(\sigma_k)_{k \in
N}$ and $(\omega_k)_{k \in N}$ defined by
$(6.2)-(6.3)$ are convergent.

Then  the solution $(\Phi_t)_{t \in R}$ of $(6.4)-(6.5)$ is defined by
$$\Phi_t((C_k)_{k \in N})=
e^{t(\sum_{n=0}^{\infty}
A_n\Big((\mathcal{F})\frac{\partial}{\partial x}\Big)^{n})}\times
(C_k)_{k \in N} \eqno(6.6)
$$
where $e^{t(\sum_{n=0}^{\infty}
A_n\Big((\mathcal{F})\frac{\partial}{\partial x}\Big)^{n})}$
denotes  an exponent of the matrix $t(\sum_{n=0}^{\infty}
A_n\Big((\mathcal{F})\frac{\partial}{\partial x}\Big)^{n})$ and it
exactly coincides with an infinite-dimensional  $(1,2,2,\dots)$
-cellular matrix $D(t)$ with cells $(D_k(t))_{k \in N}$ for which
$D_0(t)=(e^{tA_0})$ and
$$D_k(t)= e^{\sigma_k t}\begin{pmatrix}
  \cos (\omega_k t) & \sin(\omega_k t)&  \\
  -\sin(\omega_k t) & \cos (\omega_k t)&
\end{pmatrix}. \eqno(6.7)$$


By Lemma 2.2 and Corollary 6.1, one can easily establish  the
validity of the following assertions.

{\bf Observation 6.1}~{\it~ Suppose that $(\Phi_t)_{t \in R}$ is the  dynamical system in $R^{\infty}$ which comes from Corollary 6.1. Then $(\Phi_t)_{t \in R}$  is :

a) stable in the sense of an ordinary $(1,2,2, \dots)$-Lebesgue
measure $\mu_{(1,2,2, \dots)}$ in $R^{\infty}$ if and only if
$A_0+2\sum_{k=1}^{\infty}\sigma_k=0$.

b) extensible  in the sense of an ordinary $(1,2,2,
\dots)$-Lebesgue measure $\mu_{(1,2,2, \dots)}$ in $R^{\infty}$ if
and only if $0<A_0+2\sum_{k=1}^{\infty}\sigma_k<+\infty$.

c) pressing  in the sense of an ordinary $(1,2,2, \dots)$-Lebesgue
measure $\mu_{(1,2,2, \dots)}$ in $R^{\infty}$ if and only if
$-\infty<A_0+2\sum_{k=1}^{\infty}\sigma_k<0$.

d) stable in the sense of a standard $(1,2,2, \dots)$-Lebesgue
measure $\nu_{(1,2,2, \dots)}$  in $R^{\infty}$ if and only if
$A_0+2\sum_{k=1}^{\infty}\sigma_k=0$ and the series
$A_0+2\sum_{k=1}^{\infty}\sigma_k$ is absolutely convergent.

e) extensible  in the sense of a standard  $(1,2,2,
\dots)$-Lebesgue measure $\nu_{(1,2,2, \dots)}$ in $R^{\infty}$ if
and only if $0<A_0+2\sum_{k=1}^{\infty}\sigma_k<+\infty$ and the
series $A_0+2\sum_{k=1}^{\infty}\sigma_k$ is absolutely
convergent.

f) pressing  in the sense of a standard $(1,2,2, \dots)$-Lebesgue
measure $\nu_{(1,2,2, \dots)}$ in $R^{\infty}$ if and only if
$-\infty<A_0+2\sum_{k=1}^{\infty}\sigma_k<0$ and the series
$A_0+2\sum_{k=1}^{\infty}\sigma_k$ is absolutely convergent.

g) totally extensible  in the sense of a standard  $(1,2,2,
\dots)$-Lebesgue measure $\nu_{(1,2,2, \dots)}$ in $R^{\infty}$ if
and only if the series $A_0+2\sum_{k=1}^{\infty}\sigma_k$ is not
absolutely convergent and $\sum_{k \in S_{-}}^{\infty}\sigma_k
>-\infty$, where $S_{-}$ denotes a set of all natural numbers for
which $\sigma_k<0$.

h) totally pressing  in the sense of a standard  $(1,2,2,
\dots)$-Lebesgue measure $\nu_{(1,2,2, \dots)}$ in $R^{\infty}$ if
and only if the series $A_0+2\sum_{k=1}^{\infty}\sigma_k$ is not
absolutely convergent and $\sum_{k \in S_{-}}^{\infty}\sigma_k
=-\infty$, where $S_{-}$ denotes a set of all natural numbers for
which $\sigma_k<0$.}

{\bf Satisfaction Problem  of Consumers  Demands for dynamical system $(\Phi_t)_{t \in R}$ in $R^{\infty}$ defined by $(6.6)$.}

 \medskip

 Let consider $(1,2,2, \cdots)$-ordinary
''Lebesgue measure''  $\mu_{(1,2,2, \cdots)}$. By Lemma 2.2 we know that $$\mu_{(1,2,2, \cdots)}(\Phi_t(X))= e^{t(A_0+2\sum_{k =0}^{\infty}\sigma_k)}\mu_{(1,2,2, \cdots)}(X). \eqno (6.8) $$

{\bf Case 1.} $A_0+2\sum_{k =0}^{\infty}\sigma_k$ is divergent. In that case we do not know whether SPCD for dynamical system $(\Phi_t)_{t \in R}$ in $R^{\infty}$ defined by $(6.6)$ has any solution.

{\bf Case 2.} $A_0+2\sum_{k =0}^{\infty}\sigma_k=-\infty$. In that case  SPCD for dynamical system $(\Phi_t)_{t \in R}$ in $R^{\infty}$ defined by $(6.6)$ has no any solution.

{\bf Case 3.} $A_0+2\sum_{k =0}^{\infty}\sigma_k=+\infty$. In that case SPCD  for dynamical system $(\Phi_t)_{t \in R}$ in $R^{\infty}$ defined by $(6.6)$  has no any solution because if the supplier take an arbitrary
measurable system  of the positive  $\mu_{(1,1, \cdots)}$ measure, then  demands of all consumers
will be always  satisfied.

{\bf Case 4.}$A_0+2\sum_{k =0}^{\infty}\sigma_k$ is convergent and $-\infty < A_0+2\sum_{k =0}^{\infty}\sigma_k< +\infty$.

By the scheme used in Case 4 of the Section 5, one can easily show  that the solution of  SPCD  for dynamical system $(\Phi_t)_{t \in R}$ in $R^{\infty}$ defined by $(6.6)$ is given by
$$
m=\sum_{k=1}^n m_ke^{-t_k(A_0+2\sum_{k =0}^{\infty}\sigma_k)}\eqno(6.9).
$$

{\bf Observation 6.2} We have showed that Satisfaction Problem  of Consumers  Demands  for dynamical system $(\Phi_t)_{t \in R}$ in $R^{\infty}$ defined by $(6.6)$ has the solution  if $(\Phi_t)_{t \in R}$ is  pressing, expansible or stable  in the sense of the measure  $\mu_{(1,2,2, \cdots)}$. When  $(\Phi_t)_{t \in R}$ is  totally pressing or totally expansible the the same problem has no any solution.
\medskip



\end{document}